\documentclass[12pt,a4paper]{amsart}

\usepackage{amsfonts,amsmath,amssymb,amsthm}
\usepackage[latin1]{inputenc}

\newcommand{\Ff}{\mathbb{F}} 

\renewcommand{\epsilon}{\varepsilon}
\renewcommand{\le}{\leqslant}
\renewcommand{\ge}{\geqslant}

\newcommand{\nbar}{{\underline{n}}}
\newcommand{\zbar}{{\underline{z}}}

{\theoremstyle{plain}
\newtheorem{theorem}{Theorem}     
\newtheorem{lemma}[theorem]{Lemma}       
\newtheorem*{theorem*}{Theorem}     
}
{\theoremstyle{remark}

}

\begin{document}

\title[Irreducible polynomials]{Generating series for irreducible polynomials over finite fields}

\author{Arnaud Bodin}
\email{Arnaud.Bodin@math.univ-lille1.fr}
\address{Laboratoire Paul Painlev\'e, Math\'ematiques, Universit\'e 
Lille 1, 59655 Villeneuve d'Ascq Cedex, France}



 
\begin{abstract}
We count the number of irreducible polynomials in several variables of a given degree over a finite field.
The results are expressed in terms of a generating series, an exact formula and
an asymptotic approximation.
We also consider the case of the multi-degree and the case of indecomposable polynomials.
\end{abstract}

\maketitle


\section{Introduction} 
\label{sec:intro}


A well-known formula that goes back to Gauss 
estimate the number $I_n$ of monic irreducible polynomials
among all the $N_n$ monic polynomials of degree $n$ in $\Ff_q[x]$:
$$\frac{I_n}{N_n} \sim \frac 1 n.$$

A proof of this fact from the combinatorial point of view 
has been formalized  by Ph.~Flajolet, X.~Gourdon and D.~Panario in \cite{FGP}.
They apply this formalism to count irreducible polynomials in one variable
over $\Ff_q$ with the help of generating series. These generating series are convergent.
The goal of this paper is to count irreducible polynomials in several variables.
While the formalism is the same, the series are now formal power series.
At a first glance it seems a major obstruction, but despite that the series are non-convergent lots of manipulations are possible:
computation of the terms of a generating series, approximation of this term, M\"obius inversion formula,...
For comments and examples in this vein we refer to the book of Ph.~Flajolet and R.~Sedgewick 
\cite{FS} (and especially Examples I.19, I.20, II.15, Note II.28 and Appendix A).

We start by extending the introduction and recall in Section \ref{sec:formalism} the formalism 
introduced in \cite{FGP}, we consider in Section \ref{sec:degree} the case of irreducible 
polynomials in several variables
of a given degree. Then in Section \ref{sec:multidegree} the irreducible polynomials are
counted with respect to the multi-degree (or vector-degree).
Finally in Section \ref{sec:indecomposable} we deal with indecomposable polynomials.
In each case we give an exact formula, a generating series,
an asymptotic development at first order and end by an example.
References to earlier works are given for each cases after the main statements.

\textbf{Acknowledgments.} I am pleased to thank Philippe Flajolet: 
this article starts after one of his question during a visit at \textsc{inria} Rocquencourt,
he also explained me how non-convergent series can be nicely handled.


\section{Formalism} 
\label{sec:formalism}

Let us recall the formalism of \cite{FGP}.
Let $\mathcal{I}$ be a family of primitive combinatorial objects and $\omega \in \mathcal{I}$,
then $\frac{1}{1-\omega} = 1 + \omega + \omega^2+ \cdots$ 
formally represents arbitrary sequences composed of $\omega$.
The product $\mathcal{N} = \prod_{\omega \in \mathcal{I}} {\frac{1}{1-\omega}}$
formally represents all multisets (sets with possible repetitions)
composed of elements of $\mathcal{I}$.

In our situation $\mathcal{I}$ will be the set of all monic irreducible polynomials
and $\mathcal{N}$ the set of all monic polynomials. The relation $\mathcal{N} = \prod_{\omega \in \mathcal{I}} {\frac{1}{1-\omega}}$
expresses that each monic polynomial admits a unique decomposition into monic irreducible polynomials.

In order to count polynomials over $\Ff_q$ we make the following
substitution $\omega \mapsto z^{\deg \omega}$, where
$z$ is a new formal variable and $\deg \omega$ is the degree of the polynomial $\omega$.
We get the generating series:
$$N(z) = \sum_{\omega \in \mathcal{N}} z^{\deg \omega} = \sum_{n\ge1} N_n z^n, \qquad I(z) = \sum_{\omega \in \mathcal{I}} z^{\deg \omega} = \sum_{n\ge1} I_n z^n,$$ 
that count the total number of monic polynomials and monic irreducible polynomials.
The decomposition of polynomials into irreducible factors gives the relation:
$$1+N(z) = \prod_{\omega \in \mathcal{I}} {\frac{1}{1-z^{\deg \omega}}} = \prod_{n\ge 1} \frac{1}{(1-z^n)^{I_n}}.$$
By taking the logarithm we get:
$$\log (1+N(z)) = \sum_{n \ge 1} -I_n \log(1-z^n) = \sum_{n \ge 1} \frac{I(z^n)}{n}.$$

For one variable polynomials we know that the number of monic polynomials is $N_n = q^n$ so
that $1+N(z) = 1+\sum_{n\ge 1} q^n z^n = \frac{1}{1-qz}$.
Hence
\begin{equation}
  \label{eq:intro}
  \log \frac{1}{1-qz} = \sum_{n \ge 1} \frac{I(z^n)}{n},
\end{equation}
by applying M\"obius inversion formula to the (convergent) series we
get $$I(z) = \sum_{n\ge 1} \frac{\mu(n)}{n} \log \frac{1}{1-qz}.$$
If we prefer, we can only consider the coefficients of the terms $z^n$
in Equation~(\ref{eq:intro}) to find
$$\frac{q^n}{n} = \sum_{k|n} \frac{I_{n/k}}{k}\qquad \text{ and } \qquad I_n = \frac{1}{n} \sum_{k|n} \mu(k)q^{n/k}$$
by applying the classical M\"obius inversion formula ($\mu$ is M\"obius function).
From this last expression, we deduce that $I_n = \frac{q^n}{n} + O\left( \frac{q^{n/2}}{n} \right)$,
that implies the first statement of the introduction.


\section{Irreducible polynomials: degree} 
\label{sec:degree}

Let $$b_n = \binom{\nu+n}{\nu} \quad \text{and} \quad N_n = \frac{q^{b_n} - q^{b_{n-1}}}{q-1}.$$
Then $N_n$ is the number of irreducible, monic, polynomials of degree (exactly) $n$ in
$\Ff_q[x_1,\ldots,x_\nu]$.
The corresponding generating series is
$$N(z)  = \sum_{n \ge 1} N_n z^n.$$

Let $I_n$ be the number irreducible monic polynomials of degree $n$ in $\Ff_q[x_1,\ldots,x_\nu]$,
we define the corresponding generating series:
$$I(z) = \sum_{n \ge 1} I_n z^n.$$

\begin{theorem}
\label{th:degree}
\ 
\begin{enumerate}
 \item $I(z) = \sum_{n \ge 1} \frac{\mu(n)}{n} \log(1+N(z^n))$,
 \item $I_n = \sum_{k|n} \frac{\mu(k)}{k} [z^{\frac n k}] \log(1+N(z))$,
 \item $I_n = N_n - N_1 \cdot N_{n-1}+ O(q^{b_{n-2}+b_2-2})$.
\end{enumerate}
\end{theorem}
In item (3), $\nu$ and $q$ are fixed while the error term $O$ is for values of $n$ that tend to infinity.
The history for item (3) starts with L.~Carlitz \cite{Ca1}
who proved that almost all polynomials in several variables are irreducible.
His work has been extended to the case of the bi-degree in \cite{Ca2} (as considered in Section \ref{sec:multidegree})
and by S.~Cohen \cite{Co1} for more variables. 
Here we get that most of polynomials are irreducible and an estimation of the error term.
A similar formula has been obtained in \cite{Bo}: $1-\frac{I_n}{N_n} \sim \frac{N_1\cdot N_{n-1}}{N_n}$
as $n$ grows to infinity.
A formula with explicit constants is given in \cite{vzG1} 
for polynomials in two variables $\left|1-\frac{I_n}{N_n} - q^{-n+1}\right| \le 2q^{-n}$,
with $n\ge 6$.
An asymptotic formula of higher order is proved in \cite{HM}:
$I_n = N_n + \alpha_1 N_{n-1} + \alpha_2 N_{n-2} + \cdots + \alpha_s N_{n-s} + O(N_{n-s-1})$,
where $\alpha_i$ are constants and the order $s$ is arbitrary large.

Our proof of (3) appears to be simple compare with these references. 
Finally, in the forthcoming \cite{GVZ}, the formula $I(z)$ is also obtained
and applied to get an approximation: $I_n = N_n - q^{b_n+\nu-1}\left(1+ 2/q + O(1/q^2)\right)$,
where $\nu$ and $n\ge5$ are fixed, while $q$ grows to infinity.

The end of this section is devoted to the proof of Theorem~\ref{th:degree}.

\begin{proof} \ \\
\noindent
\textbf{Generating series.}
We follow the formalism of the introduction:
$$1 + N(z) = \prod_{n \ge 1}\frac{1}{(1-z^n)^{I_n}}.$$
Hence: 
\begin{equation}
\label{eq:log}
 L(z) := \log (1+N(z)) = \sum_{n \ge 1} -I_n \log (1-z^n) = \sum_{n \ge 1} \frac{I(z^n)}{n}.
\end{equation}
By M\"obius inversion formula applied to the series (see comments after the proof) we get:
$$I(z) = \sum_{n \ge 1} \frac{\mu(n)}{n} L(z^n).$$
Another point of view is to make computations term by term: we take the coefficients of $z^n$ in Equation (\ref{eq:log})
and get
$$[z^n] L(z) = \sum_{k|n} \frac{I_{n/k}}{k}.$$
Then we apply the usual M\"obius inversion formula to $n [z^n] L(z)$ and find
\begin{equation}
\label{eq:In}
I_n = \sum_{k|n} \frac{\mu(k)}{k} [z^{\frac n k}] L(z).
\end{equation}

\noindent
\textbf{Asymptotic behaviour.}
We start with the computation of an approximation of $N(z)^2$.
\begin{align*}
[z^n]N(z)^2 
  &= [z^n]\left( N_1 z + N_2 z^2 + \cdots \right)^2 \\
  &= 2 N_{n-1}\cdot N_1 + O(q^{b_{n-2}+b_2-2}).  \\
\end{align*}
And for powers $k\ge 3$ we have:
$$[z^n]N(z)^k = O(q^{b_{n-2}+b_2-2}).$$ 
Then we compute the asymptotic behaviour
of $[z^n] L(z)$.
\begin{align*}
[z^n] L(z)
   &= [z^n] \log(1+N(z)) \\
   &= [z^n] N(z) - \frac 12 [z^n] N(z)^2 + \cdots \\
   &= N_n - \frac{1}{2} \left(2 \cdot N_1  \cdot N_{n-1} + O(q^{b_{n-2}+b_2-2}) \right) + O(q^{b_{n-2}+b_2-2}) \\
   &= N_n - N_1 \cdot N_{n-1} + O(q^{b_{n-2}+b_2-2}). \\
\end{align*}

Now we get:
\begin{align*}
I_n 
  &= \sum_{k|n} \frac{\mu(k)}{k} [z^{\frac n k}] L(z) \\
  &= \frac{\mu(1)}{1}[z^n] L(z) + \sum_{k|n, k>1} \frac{\mu(k)}{k} [z^{\frac n k}] L(z) \\
  &= N_n - N_1 \cdot N_{n-1} + O(q^{b_{n-2}+b_2-2}) + O(q^{b_{n-2}}) \\
  &= N_n - N_1 \cdot N_{n-1} + O(q^{b_{n-2}+b_2-2}).\\
\end{align*}
\end{proof}

\subsection*{Example}

For example, in two variables ($\nu = 2$), we can compute the exact expression of $I_{100}$.
We compute the first $100$ coefficients of the expansion of $L(z) = \log (1+N(z))$ into powers of $z$
(by derivation an expansion of $L(z)$ can be obtained from the inverse of $1+N(z)$).
The exact value of $I_{100}$ is then computed from formula~(2) of Theorem~\ref{th:degree}:
$$I_{100} =  \frac{1}{q-1}(q^{5151} - q^{5052}-q^{5051}-q^{5050} -q^{4955}+q^{4953}+2q^{4952}+ \cdots)$$
which has a total of $4385$ monomials!
When we specialize the approximation of Theorem~\ref{th:degree} to $n=100$ we find 
that $I_{100}$ is about $\frac{1}{q-1}\left(q^{5151}-q^{5052} -q^{5051}-q^{5050}\right)$ with an error term
of magnitude $q^{4954}$.
Finally if we set $q=2$ in the exact expression, we get
$$I_{100} = 4031880625288 \ldots 8282220076800,$$
a number with $1551$ digits.

\subsection*{Inversion formula}

We shall now sketch the proof of inversion formula for formal power series.
Suppose that $f,g \in K[[z]]$ with $f(0)=0$ and
$$g(z) = \sum_{k\ge 1} {f(z^k)} = f(z)+f(z^2)+f(z^3)+\cdots$$
then
$$f(z) = \sum_{n\ge 1} {\mu(n)g(z^n)} = \mu(1)g(z)+\mu(2)g(z^2)+\mu(3)g(z^3)+\cdots$$
where $\mu(n)$ is the ordinary M\"obius function.
The proof is in fact the same as the classic one:
$$\sum_{n\ge 1} {\mu(n)g(z^n)} = \sum_{n\ge 1} {\sum_{k\ge 1} \mu(n)  {f(z^{k\cdot n})}}.$$
In this expression the coefficient of $f(z^i)$ involves only a finite number of terms and is $\sum_{j|i} \mu(j)=1$ so that 
the former sum becomes:
$$\sum_{i\ge 1} f(z^i) \sum_{j|i} \mu(j) = \sum_{i\ge 1} f(z^i).$$

Further readings for comments on computations with formal power series can be found in \cite[Appendix A]{FS}.
For convergent series, see \cite[p.~98]{Ne} and for a general version of inversion theorem, see \cite{Ro}.


\section{Irreducible polynomials: multi-degree} 
\label{sec:multidegree}

Let $N_\nbar$ be the number of irreducible, monic, polynomials  of multi-degree $\nbar = (n_1,\ldots,n_\nu)$ in
$\Ff_q[x_1,\ldots,x_\nu]$.
By \cite{Co1} or \cite{HM} we have
$$N_\nbar = \frac{1}{q-1} \sum_{(\delta_1,\ldots,\delta_\nu) \in \{0,1\}^\nu} (-1)^{\nu+\delta_1+\cdots+\delta_\nu}q^{(n_1+\delta_1)\cdots(n_\nu+\delta_\nu)}.$$

The corresponding generating series is
$$N(z_1,\ldots,z_\nu)= \sum_{(n_1,\ldots,n_\nu) \neq (0,\ldots,0)} N_{n_1,\ldots,n_\nu} \cdot z_1^{n_1}\cdots z_\nu^{n_\nu},$$
or for short
$$N(\zbar)  = \sum_{\nbar \neq 0} N_\nbar \cdot \zbar^\nbar.$$

Let $I_\nbar$ be the number of irreducible monic polynomials of degree $\nbar$ in $\Ff_q[x_1,\ldots,x_\nu]$,
we define the corresponding generating series:
$$I(\zbar) = \sum_{\nbar \neq 0} I_\nbar \cdot \zbar^\nbar.$$

\begin{theorem}
\label{th:multidegree}
\ 
\begin{enumerate}
 \item $I(\zbar) = \sum_{k \ge 1}\frac{\mu(k)}{k} \log(1+N(\zbar^k))$,
 \item 
$I_\nbar = \sum_{k| \gcd(\nbar) } \frac{\mu(k)}{k} [\zbar^{\frac \nbar k}] \log(1+N(\zbar))$,
 \item $I_\nbar = N_\nbar - N_{1,0,\ldots,0} \cdot N_{n_1-1,n_2,\ldots,n_\nu} + O(N_{0,1,0,\ldots,0} \cdot N_{n_1,n_2-1,n_3,\ldots,n_\nu}) + O(N_{2,0,\ldots,0}\cdot N_{n_1-2,n_2,\ldots,n_\nu})$ if $n_1\ge n_2 \ge \cdots \ge n_\nu$.
\end{enumerate}
\end{theorem}
For (3) X.-D.~Hou and G.~Mullen in \cite{HM} have already obtained a similar formula. 
The scheme of the proof is the same as in the previous section.

\begin{proof} \  \\
\noindent
\textbf{Generating series.}
The proof is an application of the formalism of the introduction,
(used in a similar way as in Section~\ref{sec:degree}, see also \cite{Ca2} ``heuristic proof'' of (2.1)).
$$1 + N(\zbar) = \prod_{\nbar \neq 0} \frac{1}{(1-\zbar^\nbar)^{I_\nbar}}.$$
Hence: 
\begin{equation}
\label{eq:logmulti}
 L(\zbar) := \log (1+N(\zbar)) = \sum_{\nbar \neq 0} -I_\nbar \log (1-\zbar^\nbar) 
= \sum_{k \ge 1} \frac{I(\zbar^k)}{k},
\end{equation}
where $\zbar^k = (z_1^k,\ldots,z_\nu^k)$.
Now take the coefficient of $\zbar^\nbar$ in Equation (\ref{eq:logmulti}),
we get
$$[\zbar^\nbar] L(\zbar) = \sum_{k|\gcd(\nbar)} \frac{I_{\nbar/k}}{k}.$$
By M\"obius formula applied to $\gcd(\nbar) [\zbar^\nbar] L(\zbar)$
we get
$$
I_\nbar = \sum_{k|\gcd(\nbar)} \frac{\mu(k)}{k} [\zbar^{\frac \nbar k}] L(\zbar).
$$
We can also directly apply M\"obius inversion formula to $L(\zbar)$:
$$I(\zbar) = \sum_{k \ge 1}\frac{\mu(k)}{k} L(\zbar^k).$$

\noindent
\textbf{Asymptotic behaviour.}
We suppose $n_1\ge n_2 \ge \cdots \ge n_\nu$.
We have:
\begin{align*}
I_\nbar 
  =& \sum_{k| \gcd(\nbar)} \frac{\mu(k)}{k} [\zbar^{\frac \nbar k}] L(\zbar) \\
  =& \frac{\mu(1)}{1}[\zbar^\nbar] L(\zbar) + \sum_{k|\gcd(\nbar), k>1} \frac{\mu(k)}{k} [\zbar^{\frac \nbar k}] L(\zbar) \\
  =& N_\nbar - \frac 12 \big( 2N_{1,0,\ldots,0} \cdot N_{n_1-1,n_2,\ldots,n_\nu}
+ 2N_{0,1,0,\ldots,0} \cdot N_{n_1,n_2-1,n_3,\ldots,n_\nu}+ \cdots   \\
  &+ 
 2N_{0,\ldots,0,1} \cdot N_{n_1,n_2,\ldots,n_{\nu-1},n_\nu -1} +  O(N_{1,1,,0,\ldots,0}\cdot N_{n_1-1,n_2-1,n_3,\ldots,n_\nu}) \big)  \\
  &+  O(N_{2,0,\ldots,0}\cdot N_{n_1-2,n_2,\ldots,n_\nu}) \\
  =& N_\nbar - N_{1,0,\ldots,0} \cdot N_{n_1-1,n_2,\ldots,n_\nu} + O(N_{0,1,0,\ldots,0} \cdot N_{n_1,n_2-1,n_3,\ldots,n_\nu}) + \\
  &+O(N_{2,0,\ldots,0}\cdot N_{n_1-2,n_2,\ldots,n_\nu}). \\
\end{align*}
Intuitively these formulas express that most of reducible polynomials are of type $P(x_1,\ldots,x_\nu)=(x_1+\alpha)\cdot Q(x_1,\ldots,x_\nu)$,
where $\deg_{x_1} Q = (\deg_{x_1} P) - 1$, $\alpha$ being a constant.

\end{proof}

\subsection*{Example}
For example, 
we get $I_{11,5} = \frac{1}{q-1}\left(q^{72} - q^{67}-q^{66}-q^{60}+\cdots\right)$,
while $N_{11,5}= \frac{1}{q-1}\left(q^{72}-q^{66}-q^{60}+q^{55}\right)$.
The approximation of Theorem~\ref{th:multidegree} specialized to $(11,5)$ gives that $I_{11,5}$ is about
$\frac{1}{q-1}\left(q^{72} - q^{67}-q^{66}\right)$ with an error term of magnitude $q^{61}$.
Finally if we fix $q=2$, we get $I_{11,5} = 4499945769704095481856$
and $N_{11,5} = 4647462613867219124224$.


\section{Indecomposable polynomials} 
\label{sec:indecomposable}

\subsection*{Definition and known facts}

A polynomial $P$ in $\Ff_q[x_1,\ldots,x_\nu]$
is \emph{decomposable} if there exist $h \in \Ff_q[t]$ with $\deg h \ge 2$ and 
$Q \in \Ff_q[x_1,\ldots,x_\nu]$ such that
$$P(x_1,\ldots,x_\nu) = h \circ Q (x_1,\ldots,x_\nu).$$
Otherwise $P$ is said to be \emph{indecomposable}.
For example $P(x,y)=x^3y^3+2xy+1$ is decomposable
(with $h(t) = t^3+2t+1$ and $Q(x,y)=xy$).
Notice that if $P-c$ is irreducible in $\Ff_q[x_1,\ldots,x_\nu]$ 
for some value $c \in \Ff_q$  then $P$ is indecomposable.

Indecomposable polynomials have a special interest 
in regards of Stein's theorem and provide families 
of irreducible polynomials (see \cite{St}, \cite{Na}, \cite{BDN}):
\begin{theorem*}[Stein's theorem]
If $P \in \Ff_q[x_1,\ldots,x_\nu]$ is indecomposable then
for all but finitely many values $c \in \overline{\Ff}_q$ (the algebraic closure of $\Ff_q$),
$P-c$ is irreducible in $\overline{\Ff}_q[x_1,\ldots,x_\nu]$.
Moreover the number of values $c$ such that 
$P-c$ is reducible is strictly less than $\deg P$.
\end{theorem*}

Let $J_n$ be the number of indecomposable polynomials of degree (exactly) $n$ in $\Ff_q[x_1,\ldots,x_\nu]$,
$\nu \ge 2$.
Let $\bar N_n$ be the number of polynomials of degree $n$ (not necessarily monic) in $\Ff_q[x_1,\ldots,x_\nu]$:
$$\bar N_n = q^{b_n}-q^{b_{n-1}}.$$

\bigskip

We define the following Dirichlet series:
$$N(s)=\sum_{n \ge 1}\frac{\bar N_n}{n^s}, \qquad
J(s)=\sum_{n \ge 1}\frac{J_n}{n^s}, \qquad
F(s)=\sum_{n \ge 1}\frac{q^{n-1}}{n^s}.$$

We also define a kind of M\"obius function $\mu(d,n)$:
$$\mu(d,n) = \sum_{\substack{d=d_1 | d_2 | d_3 \cdots  | d_k | d_{k+1}=n,\\ d_1<d_2<\cdots<d_k <n}} (-1)^k q^{\frac{d_2}{d_1}+\cdots+\frac{d_{k+1}}{d_k}-k},$$
where the sum is over all the iterated divisors of $n$ ($k$ not being fixed).
We define $\mu(n,n)=1$.
For a fixed $n$, we define $1<\ell<\ell'$ to be the first two divisors of $n$.

\begin{theorem}
\label{th:indecomposable}
\ 
\begin{enumerate}
 \item $J(s) = \frac{N(s)}{F(s)}$,
 \item $J_n = \sum_{d|n, 1 \le d \le n} \mu(d,n) \bar N_d$,
 \item 
If $n$ is the product of at least three prime numbers:
$$J_n = \bar N_n - q^{\ell-1}\bar N_{\frac{n}{\ell}} + O\left( q^{\ell+\ell'-2} \bar N_{\frac{n}{\ell'}}\right).$$
\end{enumerate}
\end{theorem}

The item (3) has already been obtained in \cite{BDN} and \cite{vzG2}, where some explicit constants are
computed. For example in two variables, and when $n$ is the product of at least three prime numbers, Theorem~5.1
of \cite{BDN} yields the same first order estimate of $J_n$:
$\left|\bar N_n -J_n -  q^{\ell-1} \bar N_{\frac{n}{\ell}} \frac{1-q^{-n-1}}{1-q^{-n/\ell - 1}}\right| \le n \frac{q^{\ell-1}}{q^{n/\ell}}\bar N_{\frac{n}{\ell}}\frac{1-q^{-n-1}}{1-q^{-n/\ell - 1}}$.

\begin{proof} \ \\
\noindent
\textbf{Generating series and inversion formula.}
The decomposition of a decomposable polynomial is unique once it has been  normalized, see \cite[Section 5]{BDN}:
\begin{lemma}
\label{lem:uni}
Let $ P = h \circ Q$ be a decomposition with $\deg h \ge 2$, 
$Q$ indecomposable, monic and its constant term equals zero. Then such a decomposition is unique.
\end{lemma}

By induction it implies the following lemma (see \cite{BDN}):
\begin{lemma}
\label{lem:rec}
 $$J_n = \bar N_n - \sum_{d | n, 1 \le d < n} q^{\frac{n}{d}-1}\times J_{d}.$$
\end{lemma}
From Lemma \ref{lem:rec} we get:
$$\bar N_n = \sum_{d|n, 1 \le d \le n} q^{\frac{n}{d}-1}\times J_{d}.$$
From the point of view of Dirichlet series it yields:
$$N(s) = J(s) \cdot F(s).$$
Then by induction we also deduce from Lemma \ref{lem:rec} that
$$J_n = \sum_{d|n, 1 \le d \le n} \mu(d,n) \bar N_d.$$

\noindent
\textbf{Asymptotic behaviour.}
From this formula we deduce:
\begin{align*}
J_n 
  &= \mu(n,n) \bar N_n + \mu\left({n}/{\ell},n\right) \bar N_{\frac{n}{\ell}} + O\left(\mu\left({n}/{\ell'},n\right) \bar N_{\frac{n}{\ell'}}\right) \\
  &= \bar N_n - q^{\ell-1}\bar N_{\frac{n}{\ell}} + O\left( (q^{\ell+\ell'-2}-q^{\ell-1}) \bar N_{\frac{n}{\ell'}}\right) \\
  &= N_n - q^{\ell-1}\bar N_{\frac{n}{\ell}} + O\left(q^{\ell+\ell'-2} \bar N_{\frac{n}{\ell'}}\right).\\
\end{align*}
\end{proof}

\subsection*{Example}

For example, in two variables, we have $\bar N_{100} = q^{5151}-q^{5050}$
and we compute $J_{100}$ by formula~(2) of Theorem~\ref{th:indecomposable}.
We get $J_{100} = q^{5151}-q^{5050} - q^{1327}+ q^{1276} -q^{354} + \cdots$.
While the approximation given in Theorem~\ref{th:indecomposable} is that
$J_{100}$ is about $q^{5151}-q^{5050} - q^{1327}+ q^{1276}$ with an 
error term of magnitude $q^{355}$.


\end{document}